\documentclass[reqno]{amsart}

\usepackage{amsmath,amssymb,amsthm,amsfonts,mathrsfs}
\usepackage{fullpage}
\usepackage{xcolor}
\usepackage[colorlinks=true,
linkcolor=blue,
anchorcolor=blue,
citecolor=red
]{hyperref}
\allowdisplaybreaks % allow equations break between two pages
\newtheorem{theorem}{Theorem}[section]
\newtheorem{lemma}[theorem]{Lemma}

\newtheorem*{conjecture*}{Conjecture}
\theoremstyle{definition}

\theoremstyle{remark}

\newtheorem*{remark*}{remark}
\usepackage{colonequals}

\author{Runbo Li}
\address{International Curriculum Center, The High School Affiliated to Renmin University of China, Beijing, China}
\email{runbo.li.carey@gmail.com}

\makeatletter
\@namedef{subjclassname@2020}{\textup{}2020 Mathematics Subject Classification}
\makeatother

\title[]{On Chen's theorem, Goldbach's conjecture and almost prime twins II}

\begin{document}
	
\begin{abstract}
Let $N$ denote a sufficiently large even integer and $x$ denote a sufficiently large integer, we define $D_{1,2}(N)$ as the number of primes $p$ that such that $N - p$ has at most 2 prime factors. In this paper, we show that $D_{1,2}(N) \geqslant 1.9728 \frac{C(N) N}{(\log N)^2}$, which is rather near to the asymptotic constant $2$ in Hardy--Littlewood conjecture for Goldbach's conjecture. We also get similar results on twin prime problem and additive representations of integers. The proof combines various techniques in sieve methods, such as weighted sieve, Chen's switching principle, new distribution levels proved by Lichtman and Pascadi, Chen's double sieve and Harman's sieve.
\end{abstract}

\maketitle

%\begin{center}
%\textit{Dedicated to Professor Chengdong Pan on the occasion of the 90th anniversary of his birth.} 
%\end{center}

\tableofcontents

\section{Introduction}
One of the most famous open problem in number theory is the Goldbach's conjecture, which states that any even integers can be written as the sum of two primes. Since the original conjecture is so hard, mathematicians try to consider the problem of writing a large even integer as a sum of a prime and a number with few prime factors. Let $N$ denote a sufficiently large even integer, $p$ denote a prime, and $P_{r}$ denote an integer with at most $r$ prime factors counted with multiplicity. We define
\begin{equation}
D_{1,2}(N) = \left|\left\{p : p \leqslant N, N-p=P_{2}\right\}\right|.
\end{equation}

In 1973 Chen \cite{Chen1973} established his remarkable Chen's theorem:
\begin{equation}
D_{1,2}(N) \geqslant 0.67 \frac{C(N) N}{(\log N)^2},
\end{equation}
where
\begin{equation}
C(N) = \prod_{\substack{p \mid N \\ p>2}} \frac{p-1}{p-2} \prod_{p>2}\left(1-\frac{1}{(p-1)^{2}}\right).
\end{equation}
Chen's constant 0.67 was improved successively to
$$
0.689,\ 0.7544,\ 0.81,\ 0.8285,\ 0.836,\ 0.867,\ 0.899
$$
by Halberstam and Richert \cite{HR74} \cite{Halberstam1975}, Chen \cite{Chen1978_1} \cite{Chen1978_2}, Cai and Lu \cite{CL2002}, Wu \cite{Wu2004}, Cai \cite{CAI867} and Wu \cite{Wu2008} respectively. Chen \cite{Chen9} announced a better constant 0.9, but this work has not been published.

In our 2024 preprint \cite{LRB1733}, we increase this constant to $1.733$, which almost doubles Wu's $0.899$. In the proof we use the distribution levels of Lichtman (see \cite{Lichtman3}, and \cite{Lichtman} for an earlier development of this kind of results) and complicated techniques in sieves. In this paper, by modifying the parameters used in \cite{LRB1733} and inserting more advanced techniques, we obtain the following sharper result.
\begin{theorem}\label{t1}
$$
D_{1,2}(N) \geqslant 1.9728 \frac{C(N) N}{(\log N)^2}.
$$
\end{theorem}
Our new constant $1.9728$ gives a $13.8\%$ improvement over our previous result $1.733$ and a $119\%$ refinement of Wu's prior record $0.899$. An important meaning of our new constant is that it is very close to the conjectured asymptotic constant $2$ for $D_{1,1}(N)$, the number of primes $p$ such that $N - p$ is also a prime.

Furthermore, for two relatively prime square-free positive integers $a,b$, let $M$ denote a sufficiently large integer that is relatively prime to both $a$ and $b$, $a,b < M^{\varepsilon}$ and let $M$ be even if $a$ and $b$ are both odd. Let $R_{a,b}(M)$ denote the number of primes $p$ such that $ap$ and $M-ap$ are both square-free, $b \mid (M-ap)$, and $\frac{M-ap}{b}=P_2$. In 1976, Ross [\cite{RossPhD}, Chapter 3] established that
\begin{equation}
R_{a, b}(M) \geqslant 0.608 \frac{C(abM) M}{a b(\log M)^{2}},
\end{equation}
and in \cite{LRB1733} the constant $0.608$ was improved successively to $1.733$ by essentially the same process. Now by using the new sieve process and methods in \cite{LRB}, we have the following sharper.
\begin{theorem}\label{t2}
$$
R_{a, b}(M) \geqslant 1.9728 \frac{C(abM) M}{a b(\log M)^{2}}.
$$
\end{theorem}

Another famous problem in number theory is the twin prime problem, which states that there are infinitely many prime pairs differ by $2$. Again, mathematicians consider the problem that there are infinitely many prime $p$ such that $p+2$ has few prime factors. For the twin prime problem, let $x$ denote a sufficiently large integer and define
\begin{equation}
\pi_{1,2}(x) = \left|\left\{p : p \leqslant x, p+2=P_{2}\right\}\right|.
\end{equation}
In 1973 Chen \cite{Chen1973} showed simultaneously that
\begin{equation}
\pi_{1,2}(x) \geqslant 0.335\frac{C_2 x}{(\log x)^2},
\end{equation}
where
\begin{equation}
C_2 = 2 \prod_{p>2}\left(1-\frac{1}{(p-1)^{2}}\right),
\end{equation}
and the constant 0.608 was improved successively to
$$
0.3445,\ 0.3772,\ 0.405,\ 0.71,\ 1.015,\ 1.05,\ 1.0974,\ 1.104,\ 1.123,\ 1.13
$$
by Halberstam \cite{Halberstam1975}, Chen \cite{Chen1978_1} \cite{Chen1978_2}, Fouvry and Grupp \cite{FouvryGrupp1986}, Liu \cite{Liuhongquan}, Wu \cite{Wu1990SurLS}, Cai \cite{YingchunCai2002}, Wu \cite{Wu2004}, Cai \cite{CAI867} and Cai \cite{Cai2008} respectively. 

In \cite{LRB1733} we increase this constant to $1.238$ by similar methods. Recently, Pascadi \cite{Pascadi2} got a powerful new distribution level for primes, which is quite helpful in improving the lower bound for $\pi_{1,2}(x)$. Using his new distribution results together with sieve inputs in \cite{LRB1733}, we get the following sharper.
\begin{theorem}\label{t3}
$$
\pi_{1,2}(x) \geqslant 1.2759 \frac{C_2 x}{(\log x)^2}.
$$
\end{theorem}

\section{New distribution levels}
In this section we put $A, B>0$, $\theta_0=0$, $\theta_1=\frac{7}{32}$ from Kim--Sarnak \cite{KimSarnak}, and we define the functions $\boldsymbol{\vartheta}_{\alpha}(t_{1})$ and $\boldsymbol{\vartheta}_{\alpha}(t_1, t_2, t_3)$ with $\alpha=0 \text{ or } 1$ similar to those in \cite{Lichtman3}, but with $\theta = \theta_{\alpha}$ here for $\boldsymbol{\vartheta}_{\alpha}$. We consider the analogous set of well--factorable vectors $\mathbf{D}_{r}^{well}$:
\begin{equation}
\mathbf{D}_{r}^{well}(D)=\left\{\left(D_{1}, \ldots, D_{r}\right): D_{1} \cdots D_{m-1} D_{m}^{2}<D\ \text{for all } m \leqslant r\right\}.
\end{equation}

We shall first state the distribution results for Theorem~\ref{t1}, which were proved in \cite{Lichtman3}. We remark that the maximum possible distribution level here is $\frac{19101}{32000} \approx 0.5969$. The first one is used when Chen's switching principle is not used, and the second one is used when Chen's switching principle is used.
\begin{lemma}\label{l21}
Let $\left(D_{1}, \ldots, D_{r}\right) \in \mathbf{D}_{r}^{\text{well}}(D)$ and write $D=N^{\boldsymbol{\vartheta}}, D_{i}=N^{t_{i}}$ for $i \leqslant r$.
If $\boldsymbol{\vartheta} \leqslant \boldsymbol{\vartheta}_{1}(t_{1})-\varepsilon$, then
\begin{equation}
\sum_{\substack{b=p_{1} \cdots p_{r} \\ D_{i}<p_{i} \leqslant D_{i}^{1+\varepsilon^{9}}}} \sum_{\substack{q = bc \leqslant D \\ c \mid P\left(p_{r}\right) \\(q, N)=1}} \widetilde{\lambda}^{\pm}(q)\left(\pi(N; q, N)-\frac{\pi(N)}{\varphi(q)}\right) \ll \frac{N}{(\log N)^{A}}. \tag{i}
\end{equation}
Moreover if $t_{1} \leqslant \frac{1-\theta_1}{4}$ and $r \geqslant 3$, then (i) holds if $\boldsymbol{\vartheta} \leqslant \boldsymbol{\vartheta}_{1}(t_{1}, t_{2}, t_{3})-\varepsilon$.

If $\boldsymbol{\vartheta} \leqslant \boldsymbol{\vartheta}_{1}(t_{1})-\varepsilon$ and $r = 2$, then
\begin{equation}
\sum_{\substack{b = p_{1} p_{2} \\ D_{1}<p_{1} \leqslant D_{1}^{1+\varepsilon^{9}} \\ D_{2}<p_{2} \leqslant D_{2}^{1+\varepsilon^{9}}}} \sum_{\substack{q = bc \leqslant D \\ c \mid P\left(N^{u}\right) \\(q, N)=1}} \widetilde{\lambda}^{\pm}(q)\left(\pi(N ; q, N)-\frac{\pi(N)}{\varphi(q)}\right) \ll \frac{N}{(\log N)^{A}}. \tag{ii}
\end{equation}
Moreover if $t_{1} \leqslant \frac{1-\theta_1}{4}$, then (ii) holds if $\boldsymbol{\vartheta} \leqslant \boldsymbol{\vartheta}_{1}(t_{1}, t_{2}, u)-\varepsilon$.

If $\boldsymbol{\vartheta} \leqslant \boldsymbol{\vartheta}_{1}(t_{1})-\varepsilon$ and $r = 1$, then
\begin{equation}
\sum_{\substack{b = p_{1} \\ D_{1}<p_{1} \leqslant D_{1}^{1+\varepsilon^{9}}}} \sum_{\substack{q = bc \leqslant D \\ c \mid P\left(N^{u}\right) \\(q, N)=1}} \widetilde{\lambda}^{\pm}(q)\left(\pi(N ; q, N)-\frac{\pi(N)}{\varphi(q)}\right) \ll \frac{N}{(\log N)^{A}}. \tag{iii}
\end{equation}
Moreover if $t_{1} \leqslant \frac{1-\theta_1}{4}$, then (iii) holds if $\boldsymbol{\vartheta} \leqslant \boldsymbol{\vartheta}_{1}(t_{1}, u, u)-\varepsilon$.

If $r=0$ and $u=\frac{1}{500}$, this simplifies as
$$
\sum_{\substack{q \leqslant N^{\frac{19101}{32000}} \\ q \mid P\left(N^{1/500}\right) \\(q, N)=1}} \widetilde{\lambda}^{\pm}(q)\left(\pi(N ; q, N)-\frac{\pi(N)}{\varphi(q)}\right) \ll \frac{N}{(\log N)^{A}}.
$$
\end{lemma}

\begin{lemma}\label{l22}
Let $\left(D_{1}, \ldots, D_{r}\right) \in \mathbf{D}_{r}^{\text{well}}(D)$ and write $D=N^{\boldsymbol{\vartheta}}, D_{i}=N^{t_{i}}$ for $i \leqslant r$.
Let $\varepsilon>0$ and real numbers $\varepsilon_{1}, \ldots, \varepsilon_{k} \geqslant \varepsilon$ such that $\sum_{i \leqslant k} \varepsilon_{i}=1$, and let $\Delta=1+(\log N)^{-B}$. If $\boldsymbol{\vartheta} \leqslant \boldsymbol{\vartheta}_{1}(t_{1})-\varepsilon$, then
\begin{equation}
\sum_{\substack{b=p_{1}^{\prime} \cdots p_{r}^{\prime} \\ D_{i}<p_{i}^{\prime} \leqslant D_{i}^{1+\varepsilon^{9}}}} \sum_{\substack{q = bc \leqslant D \\ c \mid P\left(p_{r}^{\prime}\right) \\(q, N)=1}} \widetilde{\lambda}^{\pm}(q)\left(\sum_{\substack{p_1 \cdots p_k \equiv N (\bmod q) \\ N^{\varepsilon_i}/\Delta < p_i \leqslant N^{\varepsilon_i}\ \forall i \leqslant k}} 1 -\frac{1}{\varphi(q)}\sum_{\substack{(p_1 \cdots p_k , N )=1 \\ N^{\varepsilon_i}/\Delta < p_i \leqslant N^{\varepsilon_i}\ \forall i \leqslant k}} 1 \right) \ll \frac{N}{(\log N)^{A}}. \tag{i}
\end{equation}
Moreover if $t_{1} \leqslant \frac{1-\theta_1}{4}$ and $r \geqslant 3$, then (i) holds if $\boldsymbol{\vartheta} \leqslant \boldsymbol{\vartheta}_{1}(t_{1}, t_{2}, t_{3})-\varepsilon$.

If $\boldsymbol{\vartheta} \leqslant \boldsymbol{\vartheta}_{1}(t_{1})-\varepsilon$ and $r = 2$, then
\begin{equation}
\sum_{\substack{b=p_{1}^{\prime} p_{2}^{\prime} \\ D_{1}<p_{1}^{\prime} \leqslant D_{1}^{1+\varepsilon^{9}} \\ D_{2}<p_{2}^{\prime} \leqslant D_{2}^{1+\varepsilon^{9}} }} \sum_{\substack{q = bc \leqslant D \\ c \mid P\left(N^{u}\right) \\(q, N)=1}} \widetilde{\lambda}^{\pm}(q)\left(\sum_{\substack{p_1 \cdots p_k \equiv N (\bmod q) \\ N^{\varepsilon_i}/\Delta < p_i \leqslant N^{\varepsilon_i}\ \forall i \leqslant k}} 1 -\frac{1}{\varphi(q)}\sum_{\substack{(p_1 \cdots p_k , N )=1 \\ N^{\varepsilon_i}/\Delta < p_i \leqslant N^{\varepsilon_i}\ \forall i \leqslant k}} 1 \right) \ll \frac{N}{(\log N)^{A}}. \tag{ii}
\end{equation}
Moreover if $t_{1} \leqslant \frac{1-\theta_1}{4}$, then (ii) holds if $\boldsymbol{\vartheta} \leqslant \boldsymbol{\vartheta}_{1}(t_{1}, t_{2}, u)-\varepsilon$.

If $\boldsymbol{\vartheta} \leqslant \boldsymbol{\vartheta}_{1}(t_{1})-\varepsilon$ and $r = 1$, then
\begin{equation}
\sum_{\substack{b=p_{1}^{\prime} \\ D_{1}<p_{1}^{\prime} \leqslant D_{1}^{1+\varepsilon^{9}}}} \sum_{\substack{q = bc \leqslant D \\ c \mid P\left(N^{u}\right) \\(q, N)=1}} \widetilde{\lambda}^{\pm}(q)\left(\sum_{\substack{p_1 \cdots p_k \equiv N (\bmod q) \\ N^{\varepsilon_i}/\Delta < p_i \leqslant N^{\varepsilon_i}\ \forall i \leqslant k}} 1 -\frac{1}{\varphi(q)}\sum_{\substack{(p_1 \cdots p_k , N )=1 \\ N^{\varepsilon_i}/\Delta < p_i \leqslant N^{\varepsilon_i}\ \forall i \leqslant k}} 1 \right) \ll \frac{N}{(\log N)^{A}}. \tag{iii}
\end{equation}
Moreover if $t_{1} \leqslant \frac{1-\theta_1}{4}$, then (iii) holds if $\boldsymbol{\vartheta} \leqslant \boldsymbol{\vartheta}_{1}(t_{1}, u, u)-\varepsilon$.

If $r=0$ and $u=\frac{1}{500}$, this simplifies as
$$
\sum_{\substack{q \leqslant N^{\frac{19101}{32000}} \\ q \mid P\left(N^{1/500}\right) \\(q, N)=1}} \widetilde{\lambda}^{\pm}(q)\left(\sum_{\substack{p_1 \cdots p_k \equiv N (\bmod q) \\ N^{\varepsilon_i}/\Delta < p_i \leqslant N^{\varepsilon_i}\ \forall i \leqslant k}} 1 -\frac{1}{\varphi(q)}\sum_{\substack{(p_1 \cdots p_k , N )=1 \\ N^{\varepsilon_i}/\Delta < p_i \leqslant N^{\varepsilon_i}\ \forall i \leqslant k}} 1 \right) \ll \frac{N}{(\log N)^{A}}.
$$
\end{lemma}

Next we shall state the distribution results for Theorem~\ref{t3}, which were proved in \cite{Pascadi2}. We remark that the maximum possible distribution level here is $\frac{2497}{4000} = 0.62425$. The first one is used when Chen's switching principle is not used, and the second one is used when Chen's switching principle is used.
\begin{lemma}\label{l23}
Let $\left(D_{1}, \ldots, D_{r}\right) \in \mathbf{D}_{r}^{\text{well}}(D)$ and write $D=x^{\boldsymbol{\vartheta}}, D_{i}=x^{t_{i}}$ for $i \leqslant r$.
If $\boldsymbol{\vartheta} \leqslant \boldsymbol{\vartheta}_{0}(t_{1})-\varepsilon$, then
\begin{equation}
\sum_{\substack{b=p_{1} \cdots p_{r} \\ D_{i}<p_{i} \leqslant D_{i}^{1+\varepsilon^{9}}}} \sum_{\substack{q = bc \leqslant D \\ c \mid P\left(p_{r}\right) \\(q, 2)=1}} \widetilde{\lambda}^{\pm}(q)\left(\pi(x; q, -2)-\frac{\pi(x)}{\varphi(q)}\right) \ll \frac{x}{(\log x)^{A}}. \tag{i}
\end{equation}
Moreover if $t_{1} \leqslant \frac{1-\theta_0}{4-3 \theta_0}$ and $r \geqslant 3$, then (i) holds if $\boldsymbol{\vartheta} \leqslant \boldsymbol{\vartheta}_{0}(t_{1}, t_{2}, t_{3})-\varepsilon$.

If $\boldsymbol{\vartheta} \leqslant \boldsymbol{\vartheta}_{0}(t_{1})-\varepsilon$ and $r = 2$, then
\begin{equation}
\sum_{\substack{b = p_{1} p_{2} \\ D_{1}<p_{1} \leqslant D_{1}^{1+\varepsilon^{9}} \\ D_{2}<p_{2} \leqslant D_{2}^{1+\varepsilon^{9}}}} \sum_{\substack{q = bc \leqslant D \\ c \mid P\left(x^{u}\right) \\(q, 2)=1}} \widetilde{\lambda}^{\pm}(q)\left(\pi(x ; q, -2)-\frac{\pi(x)}{\varphi(q)}\right) \ll \frac{x}{(\log x)^{A}}. \tag{ii}
\end{equation}
Moreover if $t_{1} \leqslant \frac{1-\theta_0}{4-3 \theta_0}$, then (ii) holds if $\boldsymbol{\vartheta} \leqslant \boldsymbol{\vartheta}_{0}(t_{1}, t_{2}, u)-\varepsilon$.

If $\boldsymbol{\vartheta} \leqslant \boldsymbol{\vartheta}_{0}(t_{1})-\varepsilon$ and $r = 1$, then
\begin{equation}
\sum_{\substack{b = p_{1} \\ D_{1}<p_{1} \leqslant D_{1}^{1+\varepsilon^{9}}}} \sum_{\substack{q = bc \leqslant D \\ c \mid P\left(x^{u}\right) \\(q, 2)=1}} \widetilde{\lambda}^{\pm}(q)\left(\pi(x ; q, -2)-\frac{\pi(x)}{\varphi(q)}\right) \ll \frac{x}{(\log x)^{A}}. \tag{iii}
\end{equation}
Moreover if $t_{1} \leqslant \frac{1-\theta_0}{4-3 \theta_0}$, then (iii) holds if $\boldsymbol{\vartheta} \leqslant \boldsymbol{\vartheta}_{0}(t_{1}, u, u)-\varepsilon$.

If $r=0$ and $u=\frac{1}{500}$, this simplifies as
$$
\sum_{\substack{q \leqslant x^{\frac{2497}{4000}} \\ q \mid P\left(x^{1/500}\right) \\(q, 2)=1}} \widetilde{\lambda}^{\pm}(q)\left(\pi(x ; q, -2)-\frac{\pi(x)}{\varphi(q)}\right) \ll \frac{x}{(\log x)^{A}}.
$$
\end{lemma}

\begin{lemma}\label{l24}
Let $\left(D_{1}, \ldots, D_{r}\right) \in \mathbf{D}_{r}^{\text{well}}(D)$ and write $D=x^{\boldsymbol{\vartheta}}, D_{i}=x^{t_{i}}$ for $i \leqslant r$.
Let $\varepsilon>0$ and real numbers $\varepsilon_{1}, \ldots, \varepsilon_{k} \geqslant \varepsilon$ such that $\sum_{i \leqslant k} \varepsilon_{i}=1$, and let $\Delta=1+(\log x)^{-B}$. If $\boldsymbol{\vartheta} \leqslant \boldsymbol{\vartheta}_{0}(t_{1})-\varepsilon$, then
\begin{equation}
\sum_{\substack{b=p_{1}^{\prime} \cdots p_{r}^{\prime} \\ D_{i}<p_{i}^{\prime} \leqslant D_{i}^{1+\varepsilon^{9}}}} \sum_{\substack{q = bc \leqslant D \\ c \mid P\left(p_{r}^{\prime}\right) \\(q, 2)=1}} \widetilde{\lambda}^{\pm}(q)\left(\sum_{\substack{p_1 \cdots p_k \equiv 2 (\bmod q) \\ x^{\varepsilon_i}/\Delta < p_i \leqslant x^{\varepsilon_i}\ \forall i \leqslant k}} 1 -\frac{1}{\varphi(q)}\sum_{\substack{(p_1 \cdots p_k , 2 )=1 \\ x^{\varepsilon_i}/\Delta < p_i \leqslant x^{\varepsilon_i}\ \forall i \leqslant k}} 1 \right) \ll \frac{x}{(\log x)^{A}}. \tag{i}
\end{equation}
Moreover if $t_{1} \leqslant \frac{1-\theta_0}{4-3 \theta_0}$ and $r \geqslant 3$, then (i) holds if $\boldsymbol{\vartheta} \leqslant \boldsymbol{\vartheta}_{0}(t_{1}, t_{2}, t_{3})-\varepsilon$.

If $\boldsymbol{\vartheta} \leqslant \boldsymbol{\vartheta}_{0}(t_{1})-\varepsilon$ and $r = 2$, then
\begin{equation}
\sum_{\substack{b=p_{1}^{\prime} p_{2}^{\prime} \\ D_{1}<p_{1}^{\prime} \leqslant D_{1}^{1+\varepsilon^{9}} \\ D_{2}<p_{2}^{\prime} \leqslant D_{2}^{1+\varepsilon^{9}} }} \sum_{\substack{q = bc \leqslant D \\ c \mid P\left(x^{u}\right) \\(q, 2)=1}} \widetilde{\lambda}^{\pm}(q)\left(\sum_{\substack{p_1 \cdots p_k \equiv 2 (\bmod q) \\ x^{\varepsilon_i}/\Delta < p_i \leqslant x^{\varepsilon_i}\ \forall i \leqslant k}} 1 -\frac{1}{\varphi(q)}\sum_{\substack{(p_1 \cdots p_k , 2 )=1 \\ x^{\varepsilon_i}/\Delta < p_i \leqslant x^{\varepsilon_i}\ \forall i \leqslant k}} 1 \right) \ll \frac{x}{(\log x)^{A}}. \tag{ii}
\end{equation}
Moreover if $t_{1} \leqslant \frac{1-\theta_0}{4-3 \theta_0}$, then (ii) holds if $\boldsymbol{\vartheta} \leqslant \boldsymbol{\vartheta}_{0}(t_{1}, t_{2}, u)-\varepsilon$.

If $\boldsymbol{\vartheta} \leqslant \boldsymbol{\vartheta}_{0}(t_{1})-\varepsilon$ and $r = 1$, then
\begin{equation}
\sum_{\substack{b=p_{1}^{\prime} \\ D_{1}<p_{1}^{\prime} \leqslant D_{1}^{1+\varepsilon^{9}}}} \sum_{\substack{q = bc \leqslant D \\ c \mid P\left(x^{u}\right) \\(q, 2)=1}} \widetilde{\lambda}^{\pm}(q)\left(\sum_{\substack{p_1 \cdots p_k \equiv 2 (\bmod q) \\ x^{\varepsilon_i}/\Delta < p_i \leqslant x^{\varepsilon_i}\ \forall i \leqslant k}} 1 -\frac{1}{\varphi(q)}\sum_{\substack{(p_1 \cdots p_k , 2 )=1 \\ x^{\varepsilon_i}/\Delta < p_i \leqslant x^{\varepsilon_i}\ \forall i \leqslant k}} 1 \right) \ll \frac{x}{(\log x)^{A}}. \tag{iii}
\end{equation}
Moreover if $t_{1} \leqslant \frac{1-\theta_0}{4-3 \theta_0}$, then (iii) holds if $\boldsymbol{\vartheta} \leqslant \boldsymbol{\vartheta}_{0}(t_{1}, u, u)-\varepsilon$.

If $r=0$ and $u=\frac{1}{500}$, this simplifies as
$$
\sum_{\substack{q \leqslant x^{\frac{2497}{4000}} \\ q \mid P\left(x^{1/500}\right) \\(q, 2)=1}} \widetilde{\lambda}^{\pm}(q)\left(\sum_{\substack{p_1 \cdots p_k \equiv 2 (\bmod q) \\ x^{\varepsilon_i}/\Delta < p_i \leqslant x^{\varepsilon_i}\ \forall i \leqslant k}} 1 -\frac{1}{\varphi(q)}\sum_{\substack{(p_1 \cdots p_k , 2 )=1 \\ x^{\varepsilon_i}/\Delta < p_i \leqslant x^{\varepsilon_i}\ \forall i \leqslant k}} 1 \right) \ll \frac{x}{(\log x)^{A}}.
$$
\end{lemma}

\section{Weighted sieve method}
Let $\mathcal{A}$ and $\mathcal{B}$ denote finite sets of positive integers, $\mathcal{P}$ denote an infinite set of primes and $z \geqslant 2$. Put
$$
\mathcal{A}=\left\{N-p : p \leqslant N\right\}, \quad \mathcal{B}=\left\{p+2 : p \leqslant x\right\},
$$

$$
\mathcal{P}=\{p : (p, 2)=1\},\quad
\mathcal{P}(q)=\{p : p \in \mathcal{P},(p, q)=1\},
$$

$$
P(z)=\prod_{\substack{p\in \mathcal{P}\\p<z}} p,\quad
\mathcal{A}_{d}=\{a : a \in \mathcal{A}, a \equiv 0(\bmod d)\},\quad
S(\mathcal{A}; \mathcal{P},z)=\sum_{\substack{a \in \mathcal{A} \\ (a, P(z))=1}} 1.
$$

\begin{lemma}\label{l31}
We have
\begin{align*}
4 D_{1, 2}(N) \geqslant&\ 3S\left(\mathcal{A};\mathcal{P}(N), N^{\frac{1}{11.49}}\right) +S\left(\mathcal{A};\mathcal{P}(N), N^{\frac{1}{6.18}}\right)\\
& - 2 \sum_{\substack{N^{\frac{1}{11.49}} \leqslant p<N^{\frac{25}{128}} \\ (p, N)=1 }} S\left(\mathcal{A}_{p};\mathcal{P}(N), N^{\frac{1}{11.49}}\right)\\
& - \sum_{\substack{N^{\frac{25}{128}} \leqslant p<N^{\frac{1}{4}} \\ (p, N)=1 }} S\left(\mathcal{A}_{p};\mathcal{P}(N), N^{\frac{1}{11.49}}\right)\\
& - \sum_{\substack{N^{\frac{1}{4}} \leqslant p<N^{\frac{57}{224}} \\ (p, N)=1 }} S\left(\mathcal{A}_{p};\mathcal{P}(N), N^{\frac{1}{11.49}}\right)\\
& - \sum_{\substack{N^{\frac{57}{224}} \leqslant p<N^{\frac{1}{3}} \\ (p, N)=1 }} S\left(\mathcal{A}_{p};\mathcal{P}(N), N^{\frac{1}{11.49}}\right)\\
& - \sum_{\substack{N^{\frac{25}{128}} \leqslant p<N^{\frac{1}{2}-\frac{3}{11.49}} \\ (p, N)=1 }} S\left(\mathcal{A}_{p};\mathcal{P}(N), N^{\frac{1}{11.49}}\right)\\
& + \sum_{\substack{N^{\frac{1}{11.49}} \leqslant p_2 <p_1 <N^{\frac{1}{6.18}} \\ (p_1 p_2, N)=1 }} S\left(\mathcal{A}_{p_1 p_2};\mathcal{P}(N), N^{\frac{1}{11.49}}\right)\\
& + \sum_{\substack{N^{\frac{1}{11.49}} \leqslant p_2 <N^{\frac{1}{6.18}} \leqslant p_1 < N^{\frac{25}{128}} \\ (p_1 p_2, N)=1 }} S\left(\mathcal{A}_{p_1 p_2};\mathcal{P}(N), N^{\frac{1}{11.49}}\right)\\
& + \sum_{\substack{N^{\frac{1}{11.49}} \leqslant p_2 <N^{\frac{1}{6.18}} < N^{\frac{25}{128}} \leqslant p_1 < N^{\frac{1}{2}-\frac{3}{11.49}} \\ (p_1 p_2, N)=1 }} S\left(\mathcal{A}_{p_1 p_2};\mathcal{P}(N), N^{\frac{1}{11.49}}\right)\\
& - 2\sum_{\substack{N^{\frac{1}{2}-\frac{3}{11.49}} \leqslant p_1 < p_2 <(\frac{N}{p_1})^{\frac{1}{2}} \\ (p_1 p_2, N)=1  } }S\left(\mathcal{A}_{p_1 p_2};\mathcal{P}(N p_1),p_2\right)\\
& - \sum_{\substack{N^{\frac{1}{11.49}} \leqslant p_1< N^{\frac{1}{3}} \leqslant p_2 <(\frac{N}{p_1})^{\frac{1}{2}} \\ (p_1 p_2, N)=1 }} S\left(\mathcal{A}_{p_1 p_2};\mathcal{P}(N p_1),p_2\right)\\
& - \sum_{\substack{N^{\frac{1}{6.18}} \leqslant p_1< N^{\frac{1}{2}-\frac{3}{11.49}} \leqslant p_2 <(\frac{N}{p_1})^{\frac{1}{2}} \\ (p_1 p_2, N)=1 }} S\left(\mathcal{A}_{p_1 p_2};\mathcal{P}(N p_1),\left(\frac{N}{p_1 p_2}\right)^{\frac{1}{2}}\right)\\
& - \sum_{\substack{N^{\frac{1}{11.49}} \leqslant p_4 < p_3 < p_2 < p_1 <N^{\frac{1}{6.18}} \\ (p_1 p_2 p_3 p_4, N)=1} }S\left(\mathcal{A}_{p_1 p_2 p_3 p_4};\mathcal{P}(N),p_3\right)\\
& - \sum_{\substack{N^{\frac{1}{11.49}} \leqslant p_1 < p_2 < p_3 <N^{\frac{1}{6.18}} \leqslant p_4 <N^{\frac{1}{2}-\frac{2}{11.49}}p_3^{-1} \\ (p_1 p_2 p_3 p_4, N)=1 } }S\left(\mathcal{A}_{p_1 p_2 p_3 p_4};\mathcal{P}(N),p_2\right)\\
& + O\left(N^{\frac{10.49}{11.49}}\right) \\
=&\ 3 S_1 + S_2 - 2 S_3 - S_4 - S_5 - S_6 - S_7 + S_8 + S_9 \\
&+S_{10} -2 S_{11} - S_{12} - S_{13} - S_{14} - S_{15} +O\left(N^{\frac{10.49}{11.49}}\right) .
\end{align*}
\end{lemma}
\begin{proof}
Taking $\kappa_1 = \frac{1}{11.49}$ and $\kappa_2 = \frac{1}{6.18}$ in [\cite{Wu2008}, Lemma 2.2], we get Lemma~\ref{l31}.
\end{proof}

\begin{lemma}\label{l32}
We have
\begin{align*}
4 \pi_{1, 2}(x) \geqslant&\ 3S\left(\mathcal{B};\mathcal{P}, x^{\frac{1}{12}}\right) +S\left(\mathcal{B};\mathcal{P}, x^{\frac{1}{7.2}}\right)\\
& + \sum_{\substack{x^{\frac{1}{12}} \leqslant p_2 <p_1 <x^{\frac{1}{7.2}} }} S\left(\mathcal{B}_{p_1 p_2};\mathcal{P}, x^{\frac{1}{12}}\right)\\
& + \sum_{\substack{x^{\frac{1}{12}} \leqslant p_2 <x^{\frac{1}{7.2}} \leqslant p_1 < x^{\frac{1}{4}} }} S\left(\mathcal{B}_{p_1 p_2};\mathcal{P}, x^{\frac{1}{12}}\right)\\
& + \sum_{\substack{x^{\frac{1}{12}} \leqslant p_2 <x^{\frac{1}{7.2}} <x^{\frac{1}{4}} \leqslant p_1 <\min(x^{\frac{2}{7}}, x^{\frac{17}{42}}p_2^{-1}) }} S\left(\mathcal{B}_{p_1 p_2};\mathcal{P}, x^{\frac{1}{12}}\right)\\
& - 2 \sum_{\substack{x^{\frac{1}{12}} \leqslant p<x^{\frac{1}{4}} }} S\left(\mathcal{B}_{p};\mathcal{P}, x^{\frac{1}{12}}\right)-2\sum_{\substack{x^{\frac{1}{4}} \leqslant p<x^{\frac{2}{7}-\varepsilon} }} S\left(\mathcal{B}_{p};\mathcal{P}, x^{\frac{1}{12}}\right)\\
& - \sum_{\substack{x^{\frac{2}{7}-\varepsilon} \leqslant p<x^{\frac{2}{7}} }} S\left(\mathcal{B}_{p};\mathcal{P}, x^{\frac{1}{12}}\right)-\sum_{\substack{x^{\frac{2}{7}-\varepsilon} \leqslant p<x^{\frac{29}{100}} }} S\left(\mathcal{B}_{p};\mathcal{P}, x^{\frac{1}{12}}\right)\\
& - \sum_{\substack{x^{\frac{29}{100}} \leqslant p<x^{\frac{1}{3}-\varepsilon} }} S\left(\mathcal{B}_{p};\mathcal{P}, x^{\frac{1}{12}}\right)-\sum_{\substack{x^{\frac{1}{3}-\varepsilon} \leqslant p<x^{\frac{1}{3}} }} S\left(\mathcal{B}_{p};\mathcal{P}, x^{\frac{1}{12}}\right)\\
& - \sum_{\substack{x^{\frac{1}{12}} \leqslant p_1<x^{\frac{1}{3}} \leqslant p_2 <(\frac{x}{p_1})^{\frac{1}{2}}}} S\left(\mathcal{B}_{p_1 p_2};\mathcal{P}(p_1),p_2\right)\\
& - \sum_{\substack{x^{\frac{1}{7.2}} \leqslant p_1<x^{\frac{2}{7}} \leqslant p_2 <(\frac{x}{p_1})^{\frac{1}{2}}}} S\left(\mathcal{B}_{p_1 p_2};\mathcal{P}(p_1),\left(\frac{x}{p_1 p_2}\right)^{\frac{1}{2}}\right)\\
& - 2 \sum_{\substack{x^{\frac{2}{7}} \leqslant p_1 < p_2 <(\frac{x}{p_1})^{\frac{1}{2}} } }S\left(\mathcal{B}_{p_1 p_2};\mathcal{P}(p_1),p_2\right)\\
& - \sum_{\substack{x^{\frac{1}{12}} \leqslant p_4 < p_3 < p_2 <p_1 <x^{\frac{1}{7.2}} } }S\left(\mathcal{B}_{p_1 p_2 p_3 p_4};\mathcal{P}(p_1),p_3\right)\\
& - \sum_{\substack{x^{\frac{1}{12}} \leqslant p_1 < p_2 < p_3 < x^{\frac{1}{7.2}} < p_4 < \min\left( x^{\frac{2}{7}}, x^{\frac{17}{42}} p_3^{-1} \right) } }S\left(\mathcal{B}_{p_1 p_2 p_3 p_4};\mathcal{P}(p_1),p_2\right)\\
&+ O\left(x^{\frac{11}{12}}\right) \\
=&\ 3 S^{\prime}_1 + S^{\prime}_2 + S^{\prime}_3 + S^{\prime}_4 + S^{\prime}_5 - 2 S^{\prime}_6 - 2 S^{\prime}_7 - S^{\prime}_8 - S^{\prime}_9 \\ 
& - S^{\prime}_{10} - S^{\prime}_{11} - S^{\prime}_{12} - S^{\prime}_{13} - 2 S^{\prime}_{14} - S^{\prime}_{15} - S^{\prime}_{16} + O\left(x^{\frac{11}{12}}\right).
\end{align*}
\end{lemma}
\begin{proof}
This is [\cite{Cai2008}, Lemma 3.2] and [\cite{LRB1733}, Lemma 3.2].
\end{proof}

\section{Proof of Theorem 1.1}
In this section, sets $\mathcal{A}$ and $\mathcal{P}$ are defined respectively. Let $\gamma$ denotes the Euler's constant, $F(s)$ and $f(s)$ are determined by the following differential-difference equation
\begin{align*}
\begin{cases}
F(s)=\frac{2 e^{\gamma}}{s}, \quad f(s)=0, \quad &0<s \leqslant 2,\\
(s F(s))^{\prime}=f(s-1), \quad(s f(s))^{\prime}=F(s-1), \quad &s \geqslant 2 ,
\end{cases}
\end{align*} 
and let $\omega(u)$ denotes the Buchstab function determined by the following differential-difference equation
\begin{align*}
\begin{cases}
\omega(u)=\frac{1}{u}, & \quad 1 \leqslant u \leqslant 2, \\
(u \omega(u))^{\prime}= \omega(u-1), & \quad u \geqslant 2 .
\end{cases}
\end{align*}

We first consider $S_1$ and $S_2$. By Buchstab's identity, we have
\begin{align}
\nonumber S_1 = S\left(\mathcal{A}; \mathcal{P}(N), N^{\frac{1}{11.49}}\right) =&\ S\left(\mathcal{A}; \mathcal{P}(N), N^{\frac{1}{500}}\right)-\sum_{\substack{N^{\frac{1}{500}} \leqslant p<N^{\frac{1}{11.49}} \\ (p, N)=1 }} S\left(\mathcal{A}_{p};\mathcal{P}(N), N^{\frac{1}{500}}\right) \\
\nonumber &+ \sum_{\substack{N^{\frac{1}{500}} \leqslant p_2 <p_1 <N^{\frac{1}{11.49}} \\ (p_1 p_2, N)=1 }} S\left(\mathcal{A}_{p_1 p_2 };\mathcal{P}(N), N^{\frac{1}{500}}\right) \\
&- \sum_{\substack{N^{\frac{1}{500}} \leqslant p_3< p_2 <p_1 <N^{\frac{1}{11.49}} \\ (p_1 p_2 p_3, N)=1 }} S\left(\mathcal{A}_{p_1 p_2 p_3 };\mathcal{P}(N), p_3\right)
\end{align}
and
\begin{align}
\nonumber S_2 = S\left(\mathcal{A}; \mathcal{P}(N), N^{\frac{1}{6.18}}\right) =&\ S\left(\mathcal{A}; \mathcal{P}(N), N^{\frac{1}{500}}\right)-\sum_{\substack{N^{\frac{1}{500}} \leqslant p<N^{\frac{1}{6.18}} \\ (p, N)=1 }} S\left(\mathcal{A}_{p};\mathcal{P}(N), N^{\frac{1}{500}}\right) \\
\nonumber &+ \sum_{\substack{N^{\frac{1}{500}} \leqslant p_2 <p_1 <N^{\frac{1}{6.18}} \\ (p_1 p_2, N)=1 }} S\left(\mathcal{A}_{p_1 p_2 };\mathcal{P}(N), N^{\frac{1}{500}}\right) \\
&- \sum_{\substack{N^{\frac{1}{500}} \leqslant p_3< p_2 <p_1 <N^{\frac{1}{6.18}} \\ (p_1 p_2 p_3, N)=1 }} S\left(\mathcal{A}_{p_1 p_2 p_3 };\mathcal{P}(N), p_3\right).
\end{align}
By Lemma~\ref{l21}, Iwaniec's linear sieve method and arguments in \cite{Lichtman}, \cite{Lichtman3} and \cite{LRB1733} we have
\begin{align}
\nonumber S_1 \geqslant&\ (1+o(1))\frac{2}{e^{\gamma}} \left(500 f\left(500\boldsymbol{\vartheta}_{\frac{1}{500}}\right)-500 \int_{\frac{1}{500}}^{\frac{1}{11.49}} \frac{F(500(\boldsymbol{\vartheta}_1(t, \frac{1}{500}, \frac{1}{500})-t))}{t}d t \right. \\
\nonumber & +500 \int_{\frac{1}{500}}^{\frac{1}{11.49}} \int_{\frac{1}{500}}^{t_1} \frac{f(500(\boldsymbol{\vartheta}_1(t_1, t_2, \frac{1}{500})-t_1 -t_2))}{t_1 t_2} d t_2 d t_1 \\
& \left. -\int_{\frac{1}{500}}^{\frac{1}{11.49}} \int_{\frac{1}{500}}^{t_1} \int_{\frac{1}{500}}^{t_2} \frac{F\left(\frac{(\boldsymbol{\vartheta}_1(t_1, t_2, t_3)-t_1 -t_2 -t_3)}{t_3}\right)}{t_1 t_2 t_3^{2}} d t_3 d t_2 d t_1 \right)\frac{C(N)N}{(\log N)^2}
\end{align}
and
\begin{align}
\nonumber S_2 \geqslant&\ (1+o(1))\frac{2}{e^{\gamma}} \left(500 f\left(500\boldsymbol{\vartheta}_{\frac{1}{500}}\right)-500 \int_{\frac{1}{500}}^{\frac{1}{6.18}} \frac{F(500(\boldsymbol{\vartheta}_1(t, \frac{1}{500}, \frac{1}{500})-t))}{t}d t \right. \\
\nonumber & +500 \int_{\frac{1}{500}}^{\frac{1}{6.18}} \int_{\frac{1}{500}}^{t_1} \frac{f(500(\boldsymbol{\vartheta}_1(t_1, t_2, \frac{1}{500})-t_1 -t_2))}{t_1 t_2} d t_2 d t_1 \\
& \left. -\int_{\frac{1}{500}}^{\frac{1}{6.18}} \int_{\frac{1}{500}}^{t_1} \int_{\frac{1}{500}}^{t_2} \frac{F\left(\frac{(\boldsymbol{\vartheta}_1(t_1, t_2, t_3)-t_1 -t_2 -t_3)}{t_3}\right)}{t_1 t_2 t_3^{2}} d t_3 d t_2 d t_1 \right)\frac{C(N)N}{(\log N)^2},
\end{align}
where $\boldsymbol{\vartheta}_{\frac{1}{500}}=\frac{19101}{32000}$. By numerical calculations we get that
\begin{equation}
S_1 \geqslant 12.902021 \frac{C(N)N}{(\log N)^2}
\end{equation}
and
\begin{equation}
S_2 \geqslant 6.533916 \frac{C(N)N}{(\log N)^2}.
\end{equation}

For $S_3$, we can either use Buchstab's identity and Lichtman's method to estimate $S_3$ with a better distribution level as in \cite{Lichtman3} or use Chen's double sieve technique as in \cite{Wu2008}. The first option leads to 
\begin{align}
\nonumber \sum_{\substack{p \\ (p, N)=1 }} S\left(\mathcal{A}_{p};\mathcal{P}(N), N^{\frac{1}{11.49}}\right) =&\ \sum_{\substack{p \\ (p, N)=1 }} S\left(\mathcal{A}_{p};\mathcal{P}(N), N^{\frac{1}{k}}\right) \\
\nonumber &- \sum_{\substack{ p_1 \\ N^{\frac{1}{k}} \leqslant p_2< N^{\frac{1}{11.49}} \\ (p_1 p_2, N)=1 }} S\left(\mathcal{A}_{p_1 p_2};\mathcal{P}(N), N^{\frac{1}{k}}\right) \\
&+ \sum_{\substack{ p_1 \\ N^{\frac{1}{k}} \leqslant p_3 < p_2< N^{\frac{1}{11.49}} \\ (p_1 p_2 p_3, N)=1 }} S\left(\mathcal{A}_{p_1 p_2 p_3};\mathcal{P}(N), p_3\right)
\end{align}
for some $k \geqslant 11.49$, while the second option creates a small saving on $S_3$ itself. We can also use Chen's double sieve on the first two sums on the right--hand side of (15) after applying Buchstab's identity. We don't know which of these options gives a smaller value, hence we take a minimum. By Lemma~\ref{l21}, Iwaniec's linear sieve method and arguments in \cite{Lichtman}, \cite{Lichtman3} and \cite{LRB1733} we have
\begin{align}
\nonumber S_{3} \leqslant&\ (1+o(1))\frac{2}{e^{\gamma}} \left( \int_{\frac{1}{11.49}}^{\frac{25}{128}} \min\left(11.49\frac{F(11.49(\boldsymbol{\vartheta}_1(t_1, \frac{1}{11.49}, \frac{1}{11.49})-t_1))}{t_1} \right. \right. \\
\nonumber &- \frac{22.98 e^{\gamma} H(11.49(\frac{1}{2}-t_1))}{(11.49(\frac{1}{2}-t_1))t_1}, \min_{11.49 \leqslant k \leqslant 500 }\left(k\frac{F(k(\boldsymbol{\vartheta}_1(t_1, \frac{1}{k}, \frac{1}{k})-t_1))}{t_1}\right. \\
\nonumber &- \frac{2k e^{\gamma} H(k(\frac{1}{2}-t_1))}{(k(\frac{1}{2}-t_1))t_1} - k \int_{\frac{1}{k}}^{\frac{1}{11.49}} \frac{f(k(\boldsymbol{\vartheta}_1(t_1, t_2, \frac{1}{k})-t_1 -t_2))}{t_1 t_2} d t_2 \\
\nonumber &- 2k e^{\gamma} \int_{\frac{1}{k}}^{\frac{1}{11.49}} \frac{h(k(\frac{1}{2}-t_1 -t_2))}{(k(\frac{1}{2}-t_1 -t_2))t_1 t_2} d t_2 \\
\nonumber & \left. \left. \left. + \int_{\frac{1}{k}}^{\frac{1}{11.49}} \int_{\frac{1}{k}}^{t_2} \frac{F\left(\frac{(\boldsymbol{\vartheta}_1(t_1, t_2, t_3)-t_1 -t_2 -t_3)}{t_3}\right)}{t_1 t_2 t_3^{2}} d t_3 d t_2 \right) \right) d t_1 \right) \frac{C(N)N}{(\log N)^2} \\
\leqslant&\ 10.436523 \frac{C(N)N}{(\log N)^2},
\end{align}
where we choose $k=12.3$ and $H(s)=H_{1/2}(s)$ and $h(s)=h_{1/2}(s)$ are defined as the same in \cite{Wu2008}. We have used the following lower bounds of $H(s)$ and $h(s)$ for $2.0 \leqslant s \leqslant 4.9$. These values can be found at Tables 1 and 2 of \cite{Wu2008}. We remark that we have $H_{\boldsymbol{\vartheta}}(s) \geqslant H_{1/2}(s)$ and $h_{\boldsymbol{\vartheta}}(s) \geqslant h_{1/2}(s)$ for $\boldsymbol{\vartheta}>\frac{1}{2}$.
\begin{align}
H(s) \geqslant
\begin{cases}
0.0223939, & \quad 2.0 < s \leqslant 2.2, \\
0.0217196, & \quad 2.2 < s \leqslant 2.3, \\
0.0202876, & \quad 2.3 < s \leqslant 2.4, \\
0.0181433, & \quad 2.4 < s \leqslant 2.5, \\
0.0158644, & \quad 2.5 < s \leqslant 2.6, \\
0.0129923, & \quad 2.6 < s \leqslant 2.7, \\
0.0100686, & \quad 2.7 < s \leqslant 2.8, \\
0.0078162, & \quad 2.8 < s \leqslant 2.9, \\
0.0072943, & \quad 2.9 < s \leqslant 3.0, \\
0.0061642, & \quad 3.0 < s \leqslant 3.1, \\
0.0052233, & \quad 3.1 < s \leqslant 3.2, \\
0.0044073, & \quad 3.2 < s \leqslant 3.3, \\
0.0036995, & \quad 3.3 < s \leqslant 3.4, \\
0.0030860, & \quad 3.4 < s \leqslant 3.5, \\
\end{cases} \quad
\begin{cases}
0.0025551, & \quad 3.5 < s \leqslant 3.6, \\
0.0020972, & \quad 3.6 < s \leqslant 3.7, \\
0.0017038, & \quad 3.7 < s \leqslant 3.8, \\
0.0013680, & \quad 3.8 < s \leqslant 3.9, \\
0.0010835, & \quad 3.9 < s \leqslant 4.0, \\
0.0008451, & \quad 4.0 < s \leqslant 4.1, \\
0.0006482, & \quad 4.1 < s \leqslant 4.2, \\
0.0004882, & \quad 4.2 < s \leqslant 4.3, \\
0.0003602, & \quad 4.3 < s \leqslant 4.4, \\
0.0002592, & \quad 4.4 < s \leqslant 4.5, \\
0.0001803, & \quad 4.5 < s \leqslant 4.6, \\
0.0001187, & \quad 4.6 < s \leqslant 4.7, \\
0.0000702, & \quad 4.7 < s \leqslant 4.8, \\
0.0000313, & \quad 4.8 < s \leqslant 4.9, \\
\end{cases}
\end{align}
\begin{align}
h(s) \geqslant
\begin{cases}
0.0232385, & \quad s = 2.0, \\
0.0211041, & \quad 2.0 < s \leqslant 2.1, \\
0.0191556, & \quad 2.1 < s \leqslant 2.2, \\
0.0173631, & \quad 2.2 < s \leqslant 2.3, \\
0.0157035, & \quad 2.3 < s \leqslant 2.4, \\
0.0141585, & \quad 2.4 < s \leqslant 2.5, \\
0.0127132, & \quad 2.5 < s \leqslant 2.6, \\
0.0113556, & \quad 2.6 < s \leqslant 2.7, \\
0.0100756, & \quad 2.7 < s \leqslant 2.8, \\
0.0088648, & \quad 2.8 < s \leqslant 2.9, \\
0.0077612, & \quad 2.9 < s \leqslant 3.0, \\
0.0066236, & \quad 3.0 < s \leqslant 3.1, \\
0.0055818, & \quad 3.1 < s \leqslant 3.2, \\
0.0046164, & \quad 3.2 < s \leqslant 3.3, \\
0.0037529, & \quad 3.3 < s \leqslant 3.4, \\
\end{cases} \quad
\begin{cases}
0.0030123, & \quad 3.4 < s \leqslant 3.5, \\
0.0023901, & \quad 3.5 < s \leqslant 3.6, \\
0.0018997, & \quad 3.6 < s \leqslant 3.7, \\
0.0015336, & \quad 3.7 < s \leqslant 3.8, \\
0.0012593, & \quad 3.8 < s \leqslant 3.9, \\
0.0010120, & \quad 3.9 < s \leqslant 4.0, \\
0.0008099, & \quad 4.0 < s \leqslant 4.1, \\
0.0006440, & \quad 4.1 < s \leqslant 4.2, \\
0.0005084, & \quad 4.2 < s \leqslant 4.3, \\
0.0003980, & \quad 4.3 < s \leqslant 4.4, \\
0.0003085, & \quad 4.4 < s \leqslant 4.5, \\
0.0002365, & \quad 4.5 < s \leqslant 4.6, \\
0.0001791, & \quad 4.6 < s \leqslant 4.7, \\
0.0001396, & \quad 4.7 < s \leqslant 4.8, \\
0.0000981, & \quad 4.8 < s \leqslant 4.9. \\
\end{cases}
\end{align}
Similarly, for $S_4$, $S_5$ and $S_7$ we have
\begin{align}
\nonumber S_{4} \leqslant&\ (1+o(1))\frac{2}{e^{\gamma}} \left( \int_{\frac{25}{128}}^{\frac{1}{4}} \min\left(11.49\frac{F(11.49(\boldsymbol{\vartheta}_1(t_1)-t_1))}{t_1} - \frac{22.98 e^{\gamma} H(11.49(\frac{1}{2}-t_1))}{(11.49(\frac{1}{2}-t_1))t_1}, \right. \right. \\
\nonumber & \min_{11.49 \leqslant k \leqslant 500 }\left(k\frac{F(k(\boldsymbol{\vartheta}_1(t_1)-t_1))}{t_1}- \frac{2k e^{\gamma} H(k(\frac{1}{2}-t_1))}{(k(\frac{1}{2}-t_1))t_1}\right. \\
\nonumber & \left. \left. \left. -\ \int_{\frac{1}{k}}^{\frac{1}{11.49}} \frac{f\left(\frac{(\boldsymbol{\vartheta}_1(t_1)-t_1 -t_2 )}{t_2}\right)}{t_1 t_2^{2}} d t_2 \right) \right) d t_1 \right) \frac{C(N)N}{(\log N)^2} \\
\leqslant&\ 3.311305 \frac{C(N)N}{(\log N)^2}, \\
\nonumber S_{5} \leqslant&\ (1+o(1))\frac{2}{e^{\gamma}} \left( \int_{\frac{1}{4}}^{\frac{57}{224}} \min\left(11.49\frac{F(11.49(\boldsymbol{\vartheta}_1(t_1)-t_1))}{t_1}, \right. \right. \\
\nonumber & \left. \left. \min_{11.49 \leqslant k \leqslant 500 }\left(k\frac{F(k(\boldsymbol{\vartheta}_1(t_1)-t_1))}{t_1}  - \int_{\frac{1}{k}}^{\frac{1}{11.49}} \frac{f\left(\frac{(\boldsymbol{\vartheta}_1(t_1)-t_1 -t_2 )}{t_2}\right)}{t_1 t_2^{2}} d t_2 \right) \right) d t_1 \right) \frac{C(N)N}{(\log N)^2} \\
\leqslant&\ 0.272301 \frac{C(N)N}{(\log N)^2}, \\
\nonumber S_{7} \leqslant&\ (1+o(1))\frac{2}{e^{\gamma}} \left( \int_{\frac{25}{128}}^{\frac{1}{2}-\frac{3}{11.49}} \min\left(11.49\frac{F(11.49(\boldsymbol{\vartheta}_1(t_1)-t_1))}{t_1} - \frac{22.98 e^{\gamma} H(11.49(\frac{1}{2}-t_1))}{(11.49(\frac{1}{2}-t_1))t_1}, \right. \right. \\
\nonumber & \min_{11.49 \leqslant k \leqslant 500 }\left(k\frac{F(k(\boldsymbol{\vartheta}_1(t_1)-t_1))}{t_1}- \frac{2k e^{\gamma} H(k(\frac{1}{2}-t_1))}{(k(\frac{1}{2}-t_1))t_1}\right. \\
\nonumber & \left. \left. \left. -\ \int_{\frac{1}{k}}^{\frac{1}{11.49}} \frac{f\left(\frac{(\boldsymbol{\vartheta}_1(t_1)-t_1 -t_2 )}{t_2}\right)}{t_1 t_2^{2}} d t_2 \right) \right) d t_1 \right) \frac{C(N)N}{(\log N)^2} \\
\leqslant&\ 2.659313 \frac{C(N)N}{(\log N)^2}.
\end{align}

By the classical linear sieve, for $S_6$ we have
\begin{align}
\nonumber S_{6} \leqslant&\ (1+o(1))\frac{2}{e^{\gamma}} \left( 11.49 \int_{\frac{57}{224}}^{\frac{1}{3}} \frac{F(11.49(\frac{1}{2}-t))}{t}d t\right) \frac{C(N)N}{(\log N)^2} \\
\leqslant&\ 5.259433 \frac{C(N)N}{(\log N)^2}.
\end{align}

For $S_8$--$S_{10}$ we can also use Chen's double sieve to gain some savings. Using similar methods as above together with [\cite{Wu2008}, Propositions 4.2 and 4.3], we have
\begin{align}
\nonumber S_8 \geqslant&\ (1+o(1))\frac{2}{e^{\gamma}} \left( 11.49 \int_{\frac{1}{11.49}}^{\frac{1}{6.18}} \int_{\frac{1}{11.49}}^{t_1} \frac{f(11.49(\boldsymbol{\vartheta}_1(t_1, t_2, \frac{1}{11.49})-t_1 -t_2))}{t_1 t_2} d t_2 d t_1 \right. \\
\nonumber & \left. + \ 11.49 \int_{\frac{1}{11.49}}^{\frac{1}{6.18}} \int_{\frac{1}{11.49}}^{t_1} \frac{2 e^{\gamma} h(11.49(\frac{1}{2}-t_1 -t_2))}{(11.49(\frac{1}{2}-t_1 -t_2)) t_1 t_2} d t_2 d t_1 \right) \frac{C(N)N}{(\log N)^2}\\
\geqslant&\ 2.421452 \frac{C(N)N}{(\log N)^2}, \\
\nonumber S_9 \geqslant&\ (1+o(1))\frac{2}{e^{\gamma}} \left( 11.49 \int_{\frac{1}{6.18}}^{\frac{25}{128}} \int_{\frac{1}{11.49}}^{\frac{1}{6.18}}  \frac{f(11.49(\boldsymbol{\vartheta}_1(t_1, t_2, \frac{1}{11.49})-t_1 -t_2))}{t_1 t_2} d t_2 d t_1 \right. \\
\nonumber & + 11.49 \int_{\frac{1}{6.18}}^{\frac{1}{2}-\frac{2}{6.18}} \int_{\frac{1}{11.49}}^{\frac{1}{6.18}} \frac{2 e^{\gamma} h(11.49(\frac{1}{2}-t_1 -t_2))}{(11.49(\frac{1}{2}-t_1 -t_2)) t_1 t_2} d t_2 d t_1 \\
\nonumber & \left. + \ 11.49 \int_{\frac{1}{2}-\frac{2}{6.18}}^{\frac{25}{128}} \int_{\frac{1}{11.49}}^{\frac{39}{256}} \frac{2 e^{\gamma} h(11.49(\frac{1}{2}-t_1 -t_2))}{(11.49(\frac{1}{2}-t_1 -t_2)) t_1 t_2} d t_2 d t_1 \right) \frac{C(N)N}{(\log N)^2}\\
\geqslant&\ 1.382532 \frac{C(N)N}{(\log N)^2}, \\
\nonumber S_{10} \geqslant&\ (1+o(1))\frac{2}{e^{\gamma}} \left( 11.49 \int_{\frac{25}{128}}^{\frac{1}{2}-\frac{3}{11.49}} \int_{\frac{1}{11.49}}^{\frac{1}{6.18}} \frac{f(11.49(\boldsymbol{\vartheta}_1(t_1, t_2, \frac{1}{11.49})-t_1 -t_2))}{t_1 t_2} d t_2 d t_1 \right. \\
\nonumber & \left. + \ 11.49 \int_{\frac{25}{128}}^{\frac{1}{2}-\frac{3}{11.49}} \int_{\frac{1}{11.49}}^{\frac{1.5}{11.49}} \frac{2 e^{\gamma} h(11.49(\frac{1}{2}-t_1 -t_2))}{(11.49(\frac{1}{2}-t_1 -t_2)) t_1 t_2} d t_2 d t_1 \right) \frac{C(N)N}{(\log N)^2}\\
\geqslant&\ 0.960457 \frac{C(N)N}{(\log N)^2}.
\end{align}

For the remaining terms, we can use Chen's switching principle together with Lemma~\ref{l22} to estimate them. Namely, for $S_{11}$ we have
\begin{equation}
S_{11}=\sum_{\substack{N^{\frac{1}{2}-\frac{3}{11.49}} \leqslant p_1 < p_2 <(\frac{N}{p_1})^{\frac{1}{2}} \\ (p_1 p_2, N)=1  } }S\left(\mathcal{A}_{p_1 p_2};\mathcal{P}(N p_1),p_2\right) = S\left(\mathcal{A}^{\prime};\mathcal{P}(N), N^{\frac{1}{2}}\right),
\end{equation}
where the set $\mathcal{A}^{\prime}$ is defined as
$$
\mathcal{A}^{\prime}=\left\{N - p_1 p_2 m : N^{\frac{1}{2}-\frac{3}{11.49}} \leqslant p_1 < p_2 < (N/ p_1)^{\frac{1}{2}},\ p^{\prime} \mid m \Rightarrow p^{\prime} > p_2 \text{ or } p^{\prime} = p_1 \right\}.
$$
We note that each $m$ above must be a prime number or a $P_2$ since $\frac{1}{2}-\frac{3}{11.49}>\frac{1}{5}$. By Buchstab's identity, we have
\begin{align}
\nonumber S_{11} = S\left(\mathcal{A}^{\prime};\mathcal{P}(N),N^{\frac{1}{2}}\right) \leqslant&\ S\left(\mathcal{A}^{\prime};\mathcal{P}(N),N^{\frac{25}{128}}\right) \\
\nonumber =&\ S\left(\mathcal{A}^{\prime}; \mathcal{P}(N), N^{\frac{1}{500}}\right)-\sum_{\substack{N^{\frac{1}{500}} \leqslant p^{\prime} <N^{\frac{25}{128}} \\ (p^{\prime}, N)=1 }} S\left(\mathcal{A}^{\prime}_{p^{\prime}};\mathcal{P}(N), N^{\frac{1}{500}}\right) \\
\nonumber &+ \sum_{\substack{N^{\frac{1}{500}} \leqslant p^{\prime}_2 <p^{\prime}_1 <N^{\frac{25}{128}} \\ (p^{\prime}_1 p^{\prime}_2, N)=1 }} S\left(\mathcal{A}^{\prime}_{p^{\prime}_1 p^{\prime}_2 };\mathcal{P}(N), N^{\frac{1}{500}}\right) \\
&- \sum_{\substack{N^{\frac{1}{500}} \leqslant p^{\prime}_3< p^{\prime}_2 <p^{\prime}_1 <N^{\frac{25}{128}} \\ (p^{\prime}_1 p^{\prime}_2 p^{\prime}_3, N)=1 }} S\left(\mathcal{A}^{\prime}_{p^{\prime}_1 p^{\prime}_2 p^{\prime}_3 };\mathcal{P}(N), p^{\prime}_3\right).
\end{align}
Then by Lemma~\ref{l22}, Iwaniec's linear sieve method and arguments in \cite{Lichtman}, \cite{Lichtman3} and \cite{LRB1733} we have
\begin{align}
\nonumber S_{11} \leqslant&\ (1+o(1))\frac{2 C(N) \left|\mathcal{A}^{\prime}\right|}{e^{\gamma} \log N} \left(500 F\left(500\boldsymbol{\vartheta}_{\frac{1}{500}}\right)-500 \int_{\frac{1}{500}}^{\frac{25}{128}} \frac{f(500(\boldsymbol{\vartheta}_1(t, \frac{1}{500}, \frac{1}{500})-t))}{t}d t \right. \\
\nonumber & +500 \int_{\frac{1}{500}}^{\frac{25}{128}} \int_{\frac{1}{500}}^{t_1} \frac{F(500(\boldsymbol{\vartheta}_1(t_1, t_2, \frac{1}{500})-t_1 -t_2))}{t_1 t_2} d t_2 d t_1 \\
\nonumber & \left. -\int_{\frac{1}{500}}^{\frac{25}{128}} \int_{\frac{1}{500}}^{t_1} \int_{\frac{1}{500}}^{t_2} \frac{f\left(\frac{(\boldsymbol{\vartheta}_1(t_1, t_2, t_3)-t_1 -t_2 -t_3)}{t_3}\right)}{t_1 t_2 t_3^{2}} d t_3 d t_2 d t_1 \right) \\
\nonumber \leqslant&\ (1+o(1))\frac{2 G_1}{e^{\gamma}} \left( \int_{\frac{1}{2}-\frac{3}{11.49}}^{\frac{1}{3}} \int_{t_1}^{\frac{1}{2}(1-t_1)} \frac{\omega \left(\frac{1-t_1-t_2}{t_2}\right)}{t_1 t_2^2} d t_2 d t_1 \right) \frac{C(N)N}{(\log N)^2}  \\
\leqslant&\ 1.30656 \frac{C(N)N}{(\log N)^2},
\end{align}
where
\begin{align}
\nonumber G_1 =&\ 500 F\left(500\boldsymbol{\vartheta}_{\frac{1}{500}}\right)-500 \int_{\frac{1}{500}}^{\frac{25}{128}} \frac{f(500(\boldsymbol{\vartheta}_1(t, \frac{1}{500}, \frac{1}{500})-t))}{t}d t  \\
\nonumber & +500 \int_{\frac{1}{500}}^{\frac{25}{128}} \int_{\frac{1}{500}}^{t_1} \frac{F(500(\boldsymbol{\vartheta}_1(t_1, t_2, \frac{1}{500})-t_1 -t_2))}{t_1 t_2} d t_2 d t_1 \\
\nonumber & -\int_{\frac{1}{500}}^{\frac{25}{128}} \int_{\frac{1}{500}}^{t_1} \int_{\frac{1}{500}}^{t_2} \frac{f\left(\frac{(\boldsymbol{\vartheta}_1(t_1, t_2, t_3)-t_1 -t_2 -t_3)}{t_3}\right)}{t_1 t_2 t_3^{2}} d t_3 d t_2 d t_1 \\
<&\ 6.06932.
\end{align}
Similarly, for $S_{12}$ and $S_{13}$ we have
\begin{align}
\nonumber S_{12} \leqslant&\ (1+o(1))\frac{2 G_1}{e^{\gamma}} \left( \int_{\frac{1}{11.49}}^{\frac{1}{3}} \int_{\frac{1}{3}}^{\frac{1}{2}(1-t_1)} \frac{\omega \left(\frac{1-t_1-t_2}{t_2}\right)}{t_1 t_2^2} d t_2 d t_1 \right) \frac{C(N)N}{(\log N)^2}  \\
\leqslant&\ 3.912436 \frac{C(N)N}{(\log N)^2}, \\
\nonumber S_{13} \leqslant&\ (1+o(1))\frac{2 G_1}{e^{\gamma}} \left( \int_{\frac{1}{6.18}}^{\frac{1}{2}-\frac{3}{11.49}} \int_{\frac{1}{2}-\frac{3}{11.49}}^{\frac{1}{2}(1-t_1)} \frac{1}{t_1 t_2 (1-t_1-t_2)} d t_2 d t_1 \right) \frac{C(N)N}{(\log N)^2}  \\
\leqslant&\ 2.835087 \frac{C(N)N}{(\log N)^2}.
\end{align}

For $S_{14}$ and $S_{15}$, we shall use a device that has been used a lot in Harman's sieve. Since $p_3 > p_4$, we have
\begin{align}
\nonumber & \sum_{\substack{N^{\frac{1}{11.49}} \leqslant p_4 < p_3 < p_2 < p_1 <N^{\frac{1}{6.18}} \\ (p_1 p_2 p_3 p_4, N)=1} }S\left(\mathcal{A}_{p_1 p_2 p_3 p_4};\mathcal{P}(N),p_3\right) \\
\leqslant&\ \sum_{\substack{N^{\frac{1}{11.49}} \leqslant p_4 < p_3 < p_2 < p_1 <N^{\frac{1}{6.18}} \\ (p_1 p_2 p_3 p_4, N)=1} }S\left(\mathcal{A}_{p_1 p_2 p_3 p_4};\mathcal{P}(N),p_4\right).
\end{align}
Here we can apply Lemma~\ref{l21} with $r=4$ to handle part of the sum on the right--hand side of (32) if $(D_1, \ldots, D_4) \in \mathbf{D}_{4}^{well}(D)$. We use the similar arguments as above to deal with other parts. Thus, we have
\begin{align}
\nonumber S_{14} \leqslant&\ (1+o(1)) \left( \int_{\frac{1}{11.49}}^{\frac{1}{6.18}} \int_{\frac{1}{11.49}}^{t_1} \int_{\frac{1}{11.49}}^{t_2} \int_{\frac{1}{11.49}}^{t_3} \left( \texttt{Boole}[(D_1, \ldots, D_4) \in \mathbf{D}_{4}^{well}(D)] \times \right. \right. \\
\nonumber & \qquad \qquad \qquad \min \left( \frac{2}{e^{\gamma}} \frac{F\left(\frac{(\boldsymbol{\vartheta}_1(t_1, t_2, t_3)-t_1 -t_2 -t_3 -t_4)}{t_4}\right)}{t_{1} t_{2} t_{3} t_{4}^2}, \frac{2 G_1}{e^{\gamma}} \frac{\omega\left(\frac{1-t_{1}-t_{2}-t_{3}-t_{4}}{t_{3}}\right)}{t_{1} t_{2} t_{3}^2 t_{4}} \right) \\
\nonumber & \left. \left. \qquad \qquad \qquad +\ \texttt{Boole}[(D_1, \ldots, D_4) \notin \mathbf{D}_{4}^{well}(D)] \frac{2 G_1}{e^{\gamma}} \frac{\omega\left(\frac{1-t_{1}-t_{2}-t_{3}-t_{4}}{t_{3}}\right)}{t_{1} t_{2} t_{3}^2 t_{4}} \right) d t_4 d t_3 d t_2 d t_1 \right) \frac{C(N)N}{(\log N)^2} \\
\leqslant&\ 0.193502 \frac{C(N)N}{(\log N)^2}.
\end{align}
Similarly, for $S_{15}$ we have
\begin{align}
\nonumber S_{15} \leqslant&\ (1+o(1)) \left( \int_{\frac{1}{11.49}}^{\frac{1}{6.18}} \int_{t_1}^{\frac{1}{6.18}} \int_{t_2}^{\frac{1}{6.18}} \int_{\frac{1}{6.18}}^{\frac{1}{2}-\frac{2}{11.49}-t_3} \left( \texttt{Boole}[(D_4, \ldots, D_1) \in \mathbf{D}_{4}^{well}(D)] \times \right. \right. \\
\nonumber & \qquad \qquad \qquad \min \left( \frac{2}{e^{\gamma}} \frac{F\left(\frac{(\boldsymbol{\vartheta}_1(t_4, t_3, t_2)-t_1 -t_2 -t_3 -t_4)}{t_1}\right)}{t_{1}^2 t_{2} t_{3} t_{4}}, \frac{2 G_1}{e^{\gamma}} \frac{\omega\left(\frac{1-t_{1}-t_{2}-t_{3}-t_{4}}{t_{2}}\right)}{t_{1} t_{2}^2 t_{3} t_{4}} \right) \\
\nonumber & \left. \left. \qquad \qquad \qquad +\ \texttt{Boole}[(D_4, \ldots, D_1) \notin \mathbf{D}_{4}^{well}(D)] \frac{2 G_1}{e^{\gamma}} \frac{\omega\left(\frac{1-t_{1}-t_{2}-t_{3}-t_{4}}{t_{2}}\right)}{t_{1} t_{2}^2 t_{3} t_{4}} \right) d t_4 d t_3 d t_2 d t_1 \right) \frac{C(N)N}{(\log N)^2} \\
\leqslant&\ 0.183611 \frac{C(N)N}{(\log N)^2}.
\end{align}

Finally, by Lemma~\ref{l31} and (9)--(34) we get
\begin{align}
\nonumber 4 D_{1,2}(N) \geqslant&\ (3 S_1 + S_2 + S_8 + S_9 + S_{10}) \\
\nonumber &- (2 S_3 + S_4 + S_5 + S_6 + S_7 + 2 S_{11} + S_{12} + S_{13} + S_{14} + S_{15}) \\
\nonumber \geqslant&\ 7.8912 \frac{C(N)N}{(\log N)^2},
\end{align}
$$
D_{1,2}(N) \geqslant 1.9728 \frac{C(N) N}{(\log N)^2}.
$$
Theorem~\ref{t1} is proved. Since the detail of the proof of Theorem~\ref{t2} is similar to those of Theorem~\ref{t1} and Theorem 1.1 in \cite{LRB} so we omit it in this paper.

\section{Proof of Theorem 1.3}
In this section, sets $\mathcal{B}$ and $\mathcal{P}$ are defined respectively. For $S^{\prime}_1$ and $S^{\prime}_2$, by Buchstab's identity, we have
\begin{align}
\nonumber S^{\prime}_1 = S\left(\mathcal{B}; \mathcal{P}, x^{\frac{1}{12}}\right) =&\ S\left(\mathcal{B}; \mathcal{P}, x^{\frac{1}{500}}\right)-\sum_{x^{\frac{1}{500}} \leqslant p<x^{\frac{1}{12}} } S\left(\mathcal{B}_{p};\mathcal{P}, x^{\frac{1}{500}}\right) \\
\nonumber &+ \sum_{x^{\frac{1}{500}} \leqslant p_2 < p_1 <x^{\frac{1}{12}} } S\left(\mathcal{B}_{p_1 p_2};\mathcal{P}, x^{\frac{1}{500}}\right) \\
&- \sum_{x^{\frac{1}{500}} \leqslant p_3 < p_2 < p_1 <x^{\frac{1}{12}} } S\left(\mathcal{B}_{p_1 p_2 p_3};\mathcal{P}, p_3 \right)
\end{align}
and
\begin{align}
\nonumber S^{\prime}_2 = S\left(\mathcal{B}; \mathcal{P}, x^{\frac{1}{7.2}}\right) =&\ S\left(\mathcal{B}; \mathcal{P}, x^{\frac{1}{500}}\right)-\sum_{x^{\frac{1}{500}} \leqslant p<x^{\frac{1}{7.2}} } S\left(\mathcal{B}_{p};\mathcal{P}, x^{\frac{1}{500}}\right) \\
\nonumber &+ \sum_{x^{\frac{1}{500}} \leqslant p_2 < p_1 <x^{\frac{1}{7.2}} } S\left(\mathcal{B}_{p_1 p_2};\mathcal{P}, x^{\frac{1}{500}}\right) \\
&- \sum_{x^{\frac{1}{500}} \leqslant p_3 < p_2 < p_1 <x^{\frac{1}{7.2}} } S\left(\mathcal{B}_{p_1 p_2 p_3};\mathcal{P}, p_3 \right).
\end{align}
By Lemma~\ref{l23}, Iwaniec's linear sieve method and arguments in \cite{Lichtman}, \cite{Lichtman3} and \cite{LRB1733} we have
\begin{align}
\nonumber S^{\prime}_1 \geqslant&\ (1+o(1))\frac{1}{e^{\gamma}} \left(500 f\left(500\boldsymbol{\vartheta}^{\prime}_{\frac{1}{500}}\right)-500 \int_{\frac{1}{500}}^{\frac{1}{12}} \frac{F(500(\boldsymbol{\vartheta}_0(t, \frac{1}{500}, \frac{1}{500})-t))}{t}d t \right. \\
\nonumber & +500 \int_{\frac{1}{500}}^{\frac{1}{12}} \int_{\frac{1}{500}}^{t_1} \frac{f(500(\boldsymbol{\vartheta}_0(t_1, t_2, \frac{1}{500})-t_1 -t_2))}{t_1 t_2} d t_2 d t_1 \\
\nonumber & \left. -\int_{\frac{1}{500}}^{\frac{1}{12}} \int_{\frac{1}{500}}^{t_1} \int_{\frac{1}{500}}^{t_2} \frac{F\left(\frac{(\boldsymbol{\vartheta}_0(t_1, t_2, t_3)-t_1 -t_2 -t_3)}{t_3}\right)}{t_1 t_2 t_3^{2}} d t_3 d t_2 d t_1 \right)\frac{C_2 x}{(\log x)^2} \\
\geqslant&\ 6.737439 \frac{C_2 x}{(\log x)^2}
\end{align}
and
\begin{align}
\nonumber S^{\prime}_2 \geqslant&\ (1+o(1))\frac{1}{e^{\gamma}} \left(500 f\left(500\boldsymbol{\vartheta}^{\prime}_{\frac{1}{500}}\right)-500 \int_{\frac{1}{500}}^{\frac{1}{7.2}} \frac{F(500(\boldsymbol{\vartheta}_0(t, \frac{1}{500}, \frac{1}{500})-t))}{t}d t \right. \\
\nonumber & +500 \int_{\frac{1}{500}}^{\frac{1}{7.2}} \int_{\frac{1}{500}}^{t_1} \frac{f(500(\boldsymbol{\vartheta}_0(t_1, t_2, \frac{1}{500})-t_1 -t_2))}{t_1 t_2} d t_2 d t_1 \\
\nonumber & \left. -\int_{\frac{1}{500}}^{\frac{1}{7.2}} \int_{\frac{1}{500}}^{t_1} \int_{\frac{1}{500}}^{t_2} \frac{F\left(\frac{(\boldsymbol{\vartheta}_0(t_1, t_2, t_3)-t_1 -t_2 -t_3)}{t_3}\right)}{t_1 t_2 t_3^{2}} d t_3 d t_2 d t_1 \right)\frac{C_2 x}{(\log x)^2} \\
\geqslant&\ 4.011646 \frac{C_2 x}{(\log x)^2},
\end{align}
where $\boldsymbol{\vartheta}^{\prime}_{\frac{1}{500}}=\frac{2497}{4000}$. For $S^{\prime}_3$--$S^{\prime}_7$, by Lemma~\ref{l23}, Iwaniec's linear sieve method and above discussion, we have
\begin{align}
\nonumber S^{\prime}_3 \geqslant&\ (1+o(1))\frac{1}{e^{\gamma}} \left( \int_{\frac{1}{12}}^{\frac{1}{7.2}} \int_{\frac{1}{12}}^{t_1} \max \left( 12 \frac{f(12(\boldsymbol{\vartheta}_0(t_1, t_2, \frac{1}{12})-t_1 -t_2))}{t_1 t_2}, \right. \right.\\
\nonumber & \max_{12 \leqslant k \leqslant 500}\left(k \frac{f(k(\boldsymbol{\vartheta}_0(t_1, t_2, \frac{1}{k})-t_1 -t_2))}{t_1 t_2} \right. \\
\nonumber & \left. \left. \left. - \int_{\frac{1}{k}}^{\frac{1}{12}} \frac{F\left(\frac{(\boldsymbol{\vartheta}_0(t_1, t_2, t_3)-t_1 -t_2 -t_3)}{t_3}\right)}{t_1 t_2 t_3^{2}} d t_3 \right) \right) d t_2 d t_1 \right) \frac{C_2 x}{(\log x)^2}\\
\geqslant&\ 0.875194 \frac{C_2 x}{(\log x)^2},\\
\nonumber S^{\prime}_4 \geqslant&\ (1+o(1))\frac{1}{e^{\gamma}} \left( \int_{\frac{1}{7.2}}^{\frac{1}{4}} \int_{\frac{1}{12}}^{\frac{1}{7.2}} \max \left( 12 \frac{f(12(\boldsymbol{\vartheta}_0(t_1, t_2, \frac{1}{12})-t_1 -t_2))}{t_1 t_2}, \right. \right.\\
\nonumber & \max_{12 \leqslant k \leqslant 500}\left(k \frac{f(k(\boldsymbol{\vartheta}_0(t_1, t_2, \frac{1}{k})-t_1 -t_2))}{t_1 t_2} \right. \\
\nonumber & \left. \left. \left. - \int_{\frac{1}{k}}^{\frac{1}{12}} \frac{F\left(\frac{(\boldsymbol{\vartheta}_0(t_1, t_2, t_3)-t_1 -t_2 -t_3)}{t_3}\right)}{t_1 t_2 t_3^{2}} d t_3 \right) \right) d t_2 d t_1 \right) \frac{C_2 x}{(\log x)^2}\\
\geqslant&\ 1.917212 \frac{C_2 x}{(\log x)^2},\\
\nonumber S^{\prime}_5 \geqslant&\ (1+o(1))\frac{1}{e^{\gamma}} \left(12 \int_{\frac{1}{12}}^{\frac{1}{7.2}} \int_{\frac{1}{4}}^{\min\left(\frac{2}{7}, \frac{17}{42}-t_1\right)} \frac{f(12(\boldsymbol{\vartheta}_0(t_2)-t_1 -t_2))}{t_1 t_2}d t_2 d t_1 \right) \frac{C_2 x}{(\log x)^2}\\
\geqslant&\ 0.282826 \frac{C_2 x}{(\log x)^2},\\
\nonumber S^{\prime}_6 \leqslant&\ (1+o(1))\frac{1}{e^{\gamma}} \left( \int_{\frac{1}{12}}^{\frac{1}{4}} \min \left( 12 \frac{F(12(\boldsymbol{\vartheta}_0(t_1, \frac{1}{12}, \frac{1}{12})-t_1))}{t_1}, \right. \right. \\
\nonumber & \min_{12 \leqslant k \leqslant 500} \left( k \frac{F(k(\boldsymbol{\vartheta}_0(t_1, \frac{1}{k}, \frac{1}{k})-t_1))}{t_1}-k \int_{\frac{1}{k}}^{\frac{1}{12}} \frac{f(k(\boldsymbol{\vartheta}_0(t_1, t_2, \frac{1}{k})-t_1 -t_2))}{t_1 t_2} d t_2 \right. \\
\nonumber & \left. \left. \left. + \int_{\frac{1}{k}}^{\frac{1}{12}} \int_{\frac{1}{k}}^{t_2} \frac{F\left(\frac{(\boldsymbol{\vartheta}_0(t_1, t_2, t_3)-t_1 -t_2 -t_3)}{t_3}\right)}{t_1 t_2 t_3^{2}} d t_3 d t_2 \right) \right) d t_1 \right) \frac{C_2 x}{(\log x)^2} \\
\leqslant&\ 7.410929 \frac{C_2 x}{(\log x)^2},\\
\nonumber S^{\prime}_7 \leqslant&\ (1+o(1))\frac{1}{e^{\gamma}} \left( \int_{\frac{1}{4}}^{\frac{2}{7}} \min \left( 12 \frac{F(12(\boldsymbol{\vartheta}_0(t_1)-t_1))}{t_1}, \right. \right. \\
\nonumber & \left. \left. \min_{12 \leqslant k \leqslant 500} \left( k \frac{F(k(\boldsymbol{\vartheta}_0(t_1)-t_1))}{t_1} - \int_{\frac{1}{k}}^{\frac{1}{12}} \frac{f\left(\frac{(\boldsymbol{\vartheta}_0(t_1)-t_1 -t_2 )}{t_2}\right)}{t_1 t_2^{2}} d t_2 \right) \right) d t_1 \right) \frac{C_2 x}{(\log x)^2} \\
\leqslant&\ 0.925271 \frac{C_2 x}{(\log x)^2}.
\end{align}

For $S^{\prime}_{12}$--$S^{\prime}_{16}$, by Chen's switching principle, Lemma~\ref{l24} and above arguments on estimating $S_{11}$--$S_{15}$ we have
\begin{align}
\nonumber S^{\prime}_{12} \leqslant&\ (1+o(1))\frac{G_2}{e^{\gamma}} \left( \int_{2}^{11} \frac{\log \left(2-\frac{3}{t+1} \right)}{t} d t \right) \frac{C_2 x}{(\log x)^2} \\
\leqslant&\ 1.960955 \frac{C_2 x}{(\log x)^2}, \\
\nonumber S^{\prime}_{13} \leqslant&\ (1+o(1))\frac{G_2}{e^{\gamma}} \left( \int_{2.5}^{6.2} \frac{\log \left(2.5-\frac{3.5}{t+1} \right)}{t} d t \right) \frac{C_2 x}{(\log x)^2} \\
\leqslant&\ 1.699112 \frac{C_2 x}{(\log x)^2}, \\
\nonumber S^{\prime}_{14} \leqslant&\ (1+o(1))\frac{G_2}{e^{\gamma}} \left( \int_{2}^{2.5} \frac{\log (t-1)}{t} d t \right) \frac{C_2 x}{(\log x)^2} \\
\leqslant&\ 0.152213 \frac{C_2 x}{(\log x)^2}, \\
\nonumber S^{\prime}_{15} \leqslant&\ (1+o(1)) \left( \int_{\frac{1}{12}}^{\frac{1}{7.2}} \int_{\frac{1}{12}}^{t_1} \int_{\frac{1}{12}}^{t_2} \int_{\frac{1}{12}}^{t_3} \left( \texttt{Boole}[(D_1, \ldots, D_4) \in \mathbf{D}_{4}^{well}(D)] \times \right. \right. \\
\nonumber & \qquad \qquad \qquad \min \left( \frac{1}{e^{\gamma}} \frac{F\left(\frac{(\boldsymbol{\vartheta}_0(t_1, t_2, t_3)-t_1 -t_2 -t_3 -t_4)}{t_4}\right)}{t_{1} t_{2} t_{3} t_{4}^2}, \frac{G_2}{e^{\gamma}} \frac{\omega\left(\frac{1-t_{1}-t_{2}-t_{3}-t_{4}}{t_{3}}\right)}{t_{1} t_{2} t_{3}^2 t_{4}} \right) \\
\nonumber & \left. \left. \qquad \qquad \qquad +\ \texttt{Boole}[(D_1, \ldots, D_4) \notin \mathbf{D}_{4}^{well}(D)] \frac{G_2}{e^{\gamma}} \frac{\omega\left(\frac{1-t_{1}-t_{2}-t_{3}-t_{4}}{t_{3}}\right)}{t_{1} t_{2} t_{3}^2 t_{4}} \right) d t_4 d t_3 d t_2 d t_1 \right) \frac{C_2 x}{(\log x)^2} \\
\leqslant&\ 0.031709 \frac{C_2 x}{(\log x)^2}, \\
\nonumber S^{\prime}_{16} \leqslant&\ (1+o(1)) \left( \int_{\frac{1}{12}}^{\frac{1}{7.2}} \int_{t_1}^{\frac{1}{7.2}} \int_{t_2}^{\frac{1}{7.2}} \int_{\frac{1}{7.2}}^{\min\left(\frac{2}{7}, \frac{17}{42}-t_3\right)} \left( \texttt{Boole}[(D_4, \ldots, D_1) \in \mathbf{D}_{4}^{well}(D)] \times \right. \right. \\
\nonumber & \qquad \qquad \qquad \min \left( \frac{1}{e^{\gamma}} \frac{F\left(\frac{(\boldsymbol{\vartheta}_0(t_4, t_3, t_2)-t_1 -t_2 -t_3 -t_4)}{t_1}\right)}{t_{1}^2 t_{2} t_{3} t_{4}}, \frac{G_2}{e^{\gamma}} \frac{\omega\left(\frac{1-t_{1}-t_{2}-t_{3}-t_{4}}{t_{2}}\right)}{t_{1} t_{2}^2 t_{3} t_{4}} \right) \\
\nonumber & \left. \left. \qquad \qquad \qquad +\ \texttt{Boole}[(D_4, \ldots, D_1) \notin \mathbf{D}_{4}^{well}(D)] \frac{G_2}{e^{\gamma}} \frac{\omega\left(\frac{1-t_{1}-t_{2}-t_{3}-t_{4}}{t_{2}}\right)}{t_{1} t_{2}^2 t_{3} t_{4}} \right) d t_4 d t_3 d t_2 d t_1 \right) \frac{C_2 x}{(\log x)^2} \\
\leqslant&\ 0.245969 \frac{C_2 x}{(\log x)^2},
\end{align}
where
\begin{align}
\nonumber G_2 =&\ 500 F\left(500\boldsymbol{\vartheta}^{\prime}_{\frac{1}{500}}\right)-500 \int_{\frac{1}{500}}^{\frac{1}{5}} \frac{f(500(\boldsymbol{\vartheta}_0(t, \frac{1}{500}, \frac{1}{500})-t))}{t}d t  \\
\nonumber & +500 \int_{\frac{1}{500}}^{\frac{1}{5}} \int_{\frac{1}{500}}^{t_1} \frac{F(500(\boldsymbol{\vartheta}_0(t_1, t_2, \frac{1}{500})-t_1 -t_2))}{t_1 t_2} d t_2 d t_1 \\
\nonumber & -\int_{\frac{1}{500}}^{\frac{1}{5}} \int_{\frac{1}{500}}^{t_1} \int_{\frac{1}{500}}^{t_2} \frac{f\left(\frac{(\boldsymbol{\vartheta}_0(t_1, t_2, t_3)-t_1 -t_2 -t_3)}{t_3}\right)}{t_1 t_2 t_3^{2}} d t_3 d t_2 d t_1 . \\
<&\ 5.81637.
\end{align}

For the remaining terms, by the arguments in \cite{Cai2008} and \cite{Wu2008}, we have
\begin{align}
S^{\prime}_8 \ll&\ \frac{\varepsilon C_2 x}{(\log x)^{2}}, \\
S^{\prime}_9 \leqslant&\ (1+o(1))\frac{12}{e^\gamma}\left(\int_{(\frac{11}{20}-\frac{29}{100})12}^{(\frac{4}{7}-\frac{2}{7})12}\frac{F(t)}{2 \times 12 -t} d t\right) \leqslant 0.111039 \frac{C_2 x}{(\log x)^{2}}, \\
S^{\prime}_{10} \leqslant&\ (1+o(1))\frac{12}{e^\gamma}\left(\int_{(\frac{11}{20}-\frac{1}{3})12}^{(\frac{11}{20}-\frac{29}{100})12}\frac{F(t)}{\frac{11}{20} \times 12 -t} d t\right) \leqslant 1.169696 \frac{C_2 x}{(\log x)^{2}}, \\
S^{\prime}_{11} \ll&\ \frac{\varepsilon C_2 x}{(\log x)^{2}}.
\end{align}

Finally, by Lemma~\ref{l32} and (35)--(53) we get
$$
\begin{aligned}
4 \pi_{1,2}(x) \geqslant&\ (3 S^{\prime}_1 + S^{\prime}_2 + S^{\prime}_3 + S^{\prime}_4 + S^{\prime}_5) \\
&-(2 S^{\prime}_6 + 2 S^{\prime}_7 + S^{\prime}_8 + S^{\prime}_9 + S^{\prime}_{10} + S^{\prime}_{11} + S^{\prime}_{12} + S^{\prime}_{13} \\
& \quad + 2 S^{\prime}_{14} + S^{\prime}_{15} + S^{\prime}_{16} + S^{\prime}_{17} + S^{\prime}_{18} + S^{\prime}_{19}) \\
\geqslant&\ 5.1036 \frac{C_2 x}{(\log x)^{2}},
\end{aligned}
$$
$$
\pi_{1,2}(x) \geqslant 1.2759 \frac{C_2 x}{(\log x)^{2}}.
$$
Theorem~\ref{t3} is proved.

\section*{Acknowledgements} 
The author would like to thank Jiamin Li for some helpful discussions.

\bibliographystyle{plain}
\bibliography{bib}
\end{document}